\documentclass{article}
\usepackage{amsfonts}

\usepackage{graphicx}
\usepackage{amsmath}

\newtheorem{theorem}{Theorem}[section]

\newtheorem{lemma}{Lemma}[section]

\newtheorem{proposition}[theorem]{Proposition}
\newtheorem{remark}{Remark}[section]

\begin{document}

\title{New explicit examples of fixed points of Poisson shot noise transforms}
\author{Aleksander M. Iksanov* and Che Soong Kim \\
Cybernetics Faculty, Kiev T.Shevchenko\\
National University, Ukraine and \\
Department of Industrial Engineering,\\
Sangji University, Wonju, Korea 220-702 \\
iksan@unicyb.kiev.ua and\\
dowoo@mail.sangji.ac.kr}
\date{Version: August 30, 2002 }
\maketitle

\begin{abstract}
We show that gamma distributions, generalized positive Linnik distributions,
S2 distributions are fixed points of Poisson shot noise transforms. The
corresponding response functions are identified via their inverse functions
except for some special cases when those can be obtained explicitly. As a
by-product, it is proven that log-convexity of the response function is not
necessary for selfdecomposability of non-negative Poisson shot noise
distribution. Some attention is given to perpetuities of a rather special
type which are closely related to our model. In particular, we study the
problem of their existence and uniqueness.

\end{abstract}

\section{Introduction.}

\strut \footnote{%
*Corresponding author} Let $\mathcal{P}^{+}$ be the set of all probability
distributions on the Borel subsets of \newline
$\mathbb{R}^{+}=[0,\infty )$ and $h:\mathbb{R}^{+}\rightarrow \mathbb{R}^{+}$
be a Borel measurable function which in what follows we call \emph{response
function}. Fix a probability space $(\Omega ,\mathcal{F},\mathcal{P})$. It
will be assumed throughout the paper that all random variables (r.v.'s)
involved are defined there, and this space is rich enough to accumulate
independent copies of some r.v.'s. Also from now on notation $\mu =\mathcal{L%
}(\xi )$ means that $\mu \in \mathcal{P}^{+}$ is a probability distribution
of r.v. $\xi =\xi (\omega )$, $\omega \in \Omega $. The last convention is
that we always take the distribution function of measure $\mu $ that is%
\textbf{\ } right-continuous. Let $\{\tau _{i}\},i=1,2,...$ be the points of
a Poisson flow with intensity $0<\lambda <\infty $, and $\xi $, $\xi _{1}$, $%
\xi _{2}$,... be non-negative independent identically distributed (i.i.d.)
r.v.'s., independent of the Poisson flow. \strut For a fixed function $h$,
let $\mathcal{P}_{h}^{+}$ be the subset of $\mathcal{P}^{+}$ consisting of
probability distributions of r.v. $\xi $ such that the series
\begin{equation}
\sum_{i=1}^{\infty }\xi _{i}h(\tau _{i})
\end{equation}
is well-defined in the weak convergence sense (and hence in probability and
almost surely). Recall that the probability distribution of the latter
random series when exists is called\emph{\ (Poisson)} \emph{shot noise
distribution }(SND, in short).

For a fixed $\lambda $ let us define \emph{a Poisson shot noise transform}
(SNT) $\mathbb{T}_{h,\lambda }:\mathcal{P}_{h}^{+}\rightarrow \mathcal{P}%
^{+} $ as follows
\begin{equation}
\ \ \mathbb{T}_{h,\lambda }(\mathcal{L}(\xi ))=\mathcal{L}\left(
\sum_{i=1}^{\infty }\xi _{i}h(\tau _{i})\right) \text{.}
\end{equation}
At this stage we would like to remark that non-negativity assumption of the
model above is not necessary in general. It is imposed here to take into
account features of the current presentation. Iksanov, Jurek (2002b)
(henceforth to be referred to as IJ(2002)) introduce the SNT for
vector-valued response functions and distributions in many dimensions.
Furthermore Iksanov, Jurek (2002a) provide conditions on $(\mathcal{L}(\xi )$%
, $h)$ which ensure the convergence of series (1) for this more general
framework.

\strut We will say that a non-degenerate at zero probability distribution $%
\mu ^{\ast }=\mathcal{L}(\xi )$ is \emph{a fixed point} of the SNT $\mathbb{T%
}_{h,\lambda }$ and/or the pair $(\lambda $,$h)$ generates or gives rise to
a fixed point $\mu ^{\ast }$ if
\begin{equation}
\mu ^{\ast }=\mathbb{T}_{h,\lambda }(\mu ^{\ast }).
\end{equation}
Formula (3) can be rewritten in terms of the Laplace-Stieltjes transform
(LST) $\varphi ^{\ast }(s)=\int_{0}^{\infty }e^{-sx}\mu ^{\ast }(dx)$ as
follows
\begin{equation}
\varphi ^{\ast }(s)=\exp \{-\lambda \int_{0}^{\infty }(1-\varphi ^{\ast
}(sh(u)))du\}\text{.}
\end{equation}
\newline
Every Poisson SND is infinitely divisible (ID), so is $\mu ^{\ast }$.
Moreover, $\mu ^{\ast }$ has zero drift and L\'{e}vy measure $M^{\ast }$
given by its tail as follows:
\begin{equation}
M^{\ast }(x,\infty )=\lambda \int_{0}^{\infty }\mu ^{\ast }(x/h(u),\infty
)du.
\end{equation}
On the other hand by differentiating (4) (it is not hard to verify that this
is possible) and by inverting the resulting expression, one gets
\begin{equation}
\omega ^{\ast }[0,x]:=\int_{0}^{x}y\mu ^{\ast }(dy)=\int_{0}^{x}\mu ^{\ast
}[0,x-y]yM^{\ast }(dy)\text{ }
\end{equation}
\newline
(compare to standard representation of positive ID distributions due to
Steutel (1970, p.86)).\newline
Furthermore, (5) reveals that $M_{\ast }$ satisfies the relation
\begin{equation}
\int_{0}^{x}yM_{\ast }(dy)=\int_{0}^{h(+0)}\omega _{\ast }[0,x/y]\nu (dy)%
\text{, }
\end{equation}
\newline
where $\nu (dx)=-\lambda xh^{\leftarrow }(dx)$, and $h^{\leftarrow }$ is a
generalized inverse of $h$ to be defined in Section 2.\newline
Just from (4)-(7) one can deduce a lot of things about $\mu ^{\ast }$. See
Section 2 for details.

The research of fixed points of the SNT (2) has been initiated in Iksanov
(2001). There in fact the following result has been proven: if $h(x)=e^{-x}$%
, $x\geq 0$ then the condition $\lambda \leq 1$ is necessary and sufficient
to guarantee an existence of fixed points $\mu ^{\ast }$ of SNT $\mathbb{T}%
_{h,\lambda }$. Furthermore, those fixed points are positive Linnik
distributions (exponential for $\lambda =1$) which are given by the tails of
distributions $\mu ^{\ast }(x,\infty )=\sum_{k=0}^{\infty }(-\beta
)^{-k}x^{\lambda k}/\Gamma (1+\lambda k),$ \ $x\geq 0,$ $\beta >0$ , where $%
\Gamma $ stands for the Euler gamma function, or via the LST's
\begin{equation}
\int_{0}^{\infty }\exp (-zx)\mu ^{\ast }(dx)=(1+\beta z^{\lambda })^{-1}.
\end{equation}
Here it is reasonable to note that 1) Lin (2001) independently proves a
closely related result in slightly different settings by using another
approach; 2) in Iksanov (2001) the distributions with the LST (8) has been
called Mittag-Leffler distributions. However, as explained in Pakes (1995,
p. 294) (see also Lin (2001)) this may cause confusion and the name
''positive Linnik'' is more correct for these distributions.

As {it} is well-known from Vervaat (1979) or Bondesson (1992), when one
studies non-negative SND, there is no loss of generality in assuming that
the response function $h$ is right-continuous and non-increasing. Under such
assumptions IJ (2002) provide a description of fixed points that correspond
to response functions $h$ with $h(+0)\leq 1$, and also directly verify that $%
h(s)=1_{[0,a)}(s)$ for some $a>0$, and $h(s)=s^{-\alpha },\alpha >1$ give
rise to no fixed points for any positive $\lambda >0$. Also Theorem 1.1(a)
from the latter reference implies that a pair $(\lambda ,h)$ with $\lambda
\int_{0}^{\infty }h(u)du>1$ does not generate fixed points.

\strut Mentioned above are the only known before response functions which
permit either to describe fixed points explicitly (that is, to point out its
LST or distribution function etc.) or to prove an absence of fixed points.
Similarly the problem of not having many explicit examples is often
mentioned in the literature on perpetuities. This is not strange. In fact,
the reader will observe (see Lemma 3.3 below) that the size-biased
distributions which correspond to fixed points of finite mean are \emph{%
perpetuities} of a very special kind. Consequently, study of fixed points in
our model and that of perpetuities are closely related. Although those have
much in common, a certain peculiarity of fixed points requires to work out
special methods to treat them. To point out a few features of fixed points
under consideration, we only mention their ID and (in most cases) absolute
continuity on $(0,\infty )$. This is certainly not a case for general
perpetuities.

\strut Somebody may ask why one needs to seek for explicit examples of fixed
points? We believe that first it is a quite interesting theoretical problem
on its own. Second it is expected that having found the way of construction
explicit examples of fixed points, one could say more about some Lebesgue
properties of fixed points. For example, which fixed points in addition to
just mentioned ID and absolute continuity are selfdecomposable (SD)
(certainly provided that the support of $h$ is the whole half-line), or
which are unimodal? Those appear to be quite intriguing problems.

\section{\strut Main results.}

\strut Our first result states that some well-known distributions do appear
as fixed points of the SNT (2). Although Proposition 2.1 does not contain an
explicit form of the corresponding response functions except for some
partial cases (one of them can be found in the proof of Proposition 2.2), no
problems occur because the only thing one should know is that those $h$'s
are right-continuous and non-increasing with $\int_{0}^{\infty }h(u)du=1$.
Let us recall that any right continuous and non-decreasing function $g$ on $%
(0,\infty )$ allows to define its generalized inverse $g^{\leftarrow }$which
is right-continuous and non-decreasing as well and given as follows $%
g^{\leftarrow }(z)=\inf \{u:g(u)<z\}$ for $z<g(0^{+})$ and $0$ otherwise. We
also preserve the above notation for ''usual'' inverse functions which are
defined for continuous and strictly monotone $g$ by the relation $%
g(g^{\leftarrow }(z))=g^{\leftarrow }(g(z))=z$.

\begin{proposition}
a) Let $\alpha $,$\beta >0$ and $\gamma \in (0,1)$. If the function $h$ is
defined via its ''usual'' inverse
\begin{equation}
h^{\leftarrow }(x)=\alpha \int_{x}^{1}z^{-1}(1-z)^{\alpha -1}dz\text{, }x\in
(0,1)
\end{equation}
then gamma distributions $\mu _{\alpha ,\beta }(dx)=\dfrac{\beta ^{-\alpha }%
}{\Gamma (\alpha )}x^{\alpha -1}e^{-x/\beta }1_{(0,\infty )}(x)dx$ and \emph{%
generalized} positive Linnik distributions $\mu _{\alpha ,\beta ,\gamma }$
given by the LST
\begin{equation}
\int_{0}^{\infty }e^{-sx}\mu _{\alpha ,\beta ,\gamma }(dx)=\dfrac{1}{%
(1+\beta s^{\gamma })^{\alpha }}
\end{equation}
are fixed points of SNT $\mathbb{T}_{h,1\text{ }}$and $\mathbb{T}%
_{h^{1/\gamma },1\text{ }}$accordingly.$\newline
$b) Let $\delta >0$, $\rho \in (0,1)$ and $h^{\leftarrow }(x)=\ln
x+2x^{-1/2}-2$,\ $x\in (0,1)$. Then $S2$ distributions $\mu _{\delta
}(dx)=d(\sum_{n=-\infty }^{\infty }(1-2\pi ^{2}n^{2}x/\delta )e^{-\pi
^{2}n^{2}x/\delta })$ and positive distributions $\mu _{\delta ,\rho }$ with
the LST
\begin{equation}
\int_{0}^{\infty }e^{-sx}\mu _{\delta ,\rho }(dx)=\left( \dfrac{\sqrt{\delta
s^{\rho }}}{\sinh \sqrt{\delta s^{\rho }}}\right) ^{2}
\end{equation}
\newline
are fixed points of SNT's $\mathbb{T}_{h,1\text{ }}$and $\mathbb{T}%
_{h^{1/\rho },1\text{ }}$accordingly.\newline
\end{proposition}

\begin{remark}
Let $\gamma _{\alpha ,\beta }$ be a gamma r.v. and $\varepsilon \in
(0,1)\cup (2,\infty )$. Unlike the gamma distribution, $\mathcal{L}(\gamma
_{\alpha ,\beta }^{\varepsilon })$ cannot be a fixed point of SNT. If $%
\varepsilon <1$ this is so, because $\mathcal{L}(\gamma _{\alpha ,\beta
}^{\varepsilon })$ is not ID. Whereas for $\varepsilon >2$ $\mathcal{L}%
(\gamma _{1,\beta }^{\varepsilon })$ together with the lognormal
distribution are primary examples of laws which are not determined by their
moments according to Krein's criterion. The same is true for $\mathcal{L}%
(\gamma _{\alpha ,\beta }^{\varepsilon })$ as shown by Pakes, Khattree
(1992). Therefore the conclusion follows from Proposition 2.3(b).
\end{remark}

\begin{remark}
\strut All distributions of Proposition 2.1 are SD. While the background
driving L\'{e}vy processes of part a) distributions are compound Poisson
(see Iksanov, Jurek (2002a) for a recent treatment of those and
definitions), this is not the case for the others. SD of $S2$ distributions
is easy to verify because as it is shown by Pitman, Yor (2001, Table 1,
p.442) their L\'{e}vy densities are of the form $k(x)/x=(\delta
\sum_{n=-\infty }^{+\infty }e^{-\delta ^{-1}\pi ^{2}n^{2}x})/x$, and hence $%
k(x)$ is decreasing on $(0,\infty )$. Now distributions given by (11) are SD
as these are laws of strictly stable subordinator evaluated at random SD ($%
S2 $) time. The observation about SD of such distributions is due to
Bondesson (1992, p.19). \newline
Clearly, 1) NOT \emph{all} fixed points generated by $h$ of unbounded
support and 2) NOT \emph{all} size-biased distributions which correspond to
fixed points are SD.
\end{remark}

\strut To formulate our second result, recall that Bondesson (1992, p.156)
proved that the sufficient condition for SD of SND (1) is log-convexity and
strict decreasingness of $h$. The next Proposition states that this is not
necessary.

\begin{proposition}
\strut There exist selfdecomposable shot noise distributions which are
generated by a response function which is not log-convex.\strut
\end{proposition}

Suppose that $\mu \in \mathcal{P}^{+}$ is of finite mean $%
m:=\int_{0}^{\infty }x\mu (dx)$. This allows to consider the so-called
size-biased distribution $\overline{\mu }(dx)=m^{-1}x\mu (dx)$. Let $%
\overline{\eta }$, $\eta $ and $A$ be independent r.v.'s with $\overline{\mu
}=\mathcal{L}(\overline{\eta })$, $\mu =\mathcal{L}(\eta )$ and $\nu =%
\mathcal{L}(A)$ which satisfy the distributional equality
\begin{equation}
\overline{\eta }\overset{d}{=}\eta +A\overline{\eta }\text{.}
\end{equation}
\newline
We now cite the problem mentioned by Pitman, Yor (2000, p.35): ''given a
distribution of $A$...whether there exists such a distribution of $\eta $''.
Recall that in the more recent literature the so defined r.v. $\overline{%
\eta }$ (as in (12)) is typically called perpetuity.

Below we answer the above question for the partial case when $\nu $ is
concentrated on $(0,b]$, $b\leq 1$. Denote by $\delta _{x}$ the delta
measure at $x\geq 0$.

\begin{proposition}
a) For any $\nu \neq \delta _{1}$ concentrated on $(0,b]$, $b\leq 1$ there
exist $\mu $'s satisfying (12). For fixed $m>0$ $\mu $ is a unique solution
to (12) such that $m=\int_{0}^{\infty }x\mu (dx)$.\newline
b) Those $\mu $'s have finite exponential moments.\newline
c) All $\mu $'s are infinitely divisible with drift $0$ and L\'{e}vy measure
$M$ whose tail is given as follows $M(x,\infty )=\int_{0}^{b}z^{-1}\mu
(xz^{-1},\infty )\nu (dz)$. Furthermore, $\mu $'s are compound Poisson
provided $x^{-1}\nu (dx)$ is integrable at the neighbourhood of zero. \
\newline
d) If for some $\varepsilon \in (0,1]$ $\int_{0}^{b}$ $x^{-\varepsilon }\nu
(dx)<\infty $ then\newline
$\mu (dx)=q\delta _{0}+(1-q)f(x)1_{(0,\infty )}(x)dx$, where $q=0$ if $%
x^{-1}\nu (dx)$ is not integrable at the neighbourhood of the origin, and $%
q\in (0,1)$ is a unique solution to the equation $\exp (-b(1-z))=z$ if $%
\int_{0}^{b}x^{-1}\nu (dx)=b$. In words, $\mu $'s have an absolutely
continuous component on $(0,\infty )$ with density $f$.\newline
e) All $\mu $'s are fixed points of SNT $\mathbb{T}_{h,1}$ with the response
function $h$ given via its generalized inverse $h^{\leftarrow }$ as follows:
$h^{\leftarrow }(x)=\int_{x}^{b}z^{-1}\nu (dz)$ which implies $%
\int_{0}^{\infty }h(z)dz=1$.
\end{proposition}

\begin{remark}
It is possible to strengthen the above Proposition in the following way. Let
us consider the measure $\sigma $ such that $\sigma (dx)=x\mu (dx)$ and
rewrite (12) in terms of distributions to obtain the well-known
representation of positive infinitely divisible distributions due to Steutel
(1970, p.86):
\begin{equation}
\sigma \lbrack 0,x]=\int_{0}^{x}\mu \lbrack 0,x-z]zM(dz)\text{,}
\end{equation}
$M$ being the L\'{e}vy measure of $\mu $ which in our case has a feature
\begin{equation}
\int_{0}^{x}zM(dz)=\int_{0}^{b}\sigma \lbrack 0,x/z]\nu (dz).
\end{equation}
As it turned out if $\int_{0}^{b}z^{-\Delta }v(dz)<\infty $ for some$\
\Delta \in (0,1)$, we need not pre-suppose that $\int_{0}^{\infty }x\mu
(dx)<\infty $. In fact, if a distribution $\mu $ satisfies (13), (14) then
it necessarily has finite first moment. Moreover, given $m>0$ $\mu $ is the
unique distribution of mean $m$ satisfying (13), (14). This is essentially
the content of Theorem 1.1(b)\ of IJ (2002), but for a special case the
proof of that assertion should be taken into account.
\end{remark}

\section{\strut The Proofs.}

\strut Four preparatory lemmas are prepared. \strut We begin with a simple
observation which can be read from (4) and hence its proof is immediate and
omitted.\ It is singled out as a Lemma only for ease of further references.

\begin{lemma}
\strut Fixed points of the SNT (2) are invariant under scale
transformations, that is, if $\mathcal{L}(\xi )$ is a fixed point of the SNT
so is $\mathcal{L}(c\xi )$, for any $c>0$. \newline
\end{lemma}

Throughout the rest of this Section we will assume that for any positive $%
\lambda $ response functions $h$'s of the SNT $\mathbb{T}_{h,\lambda }$ are
subject to \textbf{CONDITION A}: they are right-continuous, non-increasing, $%
h(+0)\leq 1$ and $h$ is not of the form $h(u)=1_{[0,a)}(u)$ for some $a>0$.%
\newline

The next Lemma is a uniqueness result concerning fixed points of the SNT. It
is contained in Theorem 1.1(b) of IJ (2002) and has been proven there by
using Contraction Principle. We would like to provide an independent,
slightly simpler proof.

\begin{lemma}
Let $h$ satisfies Condition A and $\lambda \int_{0}^{\infty }h(z)dz=1$. Then
$\mathbb{T}_{h,\lambda }$ has fixed points of finite mean. Given $m\in
(0,\infty )$ there exist a unique fixed point $\mu ^{\ast }$ of $\mathbb{T}%
_{h,\lambda }$ with $m=\int_{0}^{\infty }x\mu ^{\ast }(dx)$.
\end{lemma}

\textbf{Proof.} For fixed $m>0$ consider the set of probability measures%
\newline
$\mathcal{P}_{h,m}^{+}=\{\rho \in \mathcal{P}_{h}^{+}:\int_{0}^{\infty
}x\rho (dx)=m\}$. Starting with $\mu _{0}=\delta _{m}$, define the sequence
\begin{equation*}
\mu _{n}:=\mathbb{T}_{h,\lambda }\mu _{n-1}:=\mathbb{T}_{h,\lambda }^{n}\mu
_{0},n=1,2,...
\end{equation*}
which is trivially well-defined on $\mathcal{P}_{h,m}^{+}$ provided $%
\int_{0}^{\infty }h(z)dz<\infty $. The corresponding LST's $\varphi
_{n}(s)=\int_{0}^{\infty }e^{-sx}\mu _{n}(dx)$, $n=0,1,...$ satisfy
equations
\begin{equation}
\varphi _{0}(s)=e^{-ms}\text{, \ \ }\varphi _{n}(s)=\exp \{-\lambda
\int_{0}^{\infty }(1-\varphi _{n-1}(sh(u)))du\}\text{, }n=1,2,...
\end{equation}
Let us verify that the weak limit of $\mu _{n}$, as $n\rightarrow \infty $,
exists and has mean $m$. As it is well-known, this will mean that $\mathbb{T}%
_{h,\lambda }$ has a unique fixed point on $\mathcal{P}_{h,m}^{+}$. \newline
In what follows we use some ideas of Durrett, Liggett (1983, the proof of
Theorem 2.7). By Jensen's inequality,
\begin{equation*}
\varphi _{1}(s)=\mathbb{E}\exp \{-s\sum_{i=1}^{\infty }\xi _{i}h(\tau
_{i})\}\geq \exp \{-s\mathbb{E}\sum_{i=1}^{\infty }\xi _{i}h(\tau
_{i})\}=\varphi _{0}(s),
\end{equation*}
that implies $\varphi _{n}(s)\geq $ $\varphi _{n-1}(s)$, $n=1,2,...$, $s\geq
0$. Thus the monotone and bounded sequence $\{\varphi _{n}\},n=1,2,...$ has
a unique limit $\varphi (s)$, say, being the LST of a probability measure $%
\mu $, say$.$ Since $\int_{0}^{\infty }h(z)dz<\infty $ then by dominated
convergence it is easily seen that $\varphi (s)$ satisfies the fixed point
equation (4) or equivalently $\mu $ is a (possibly degenerate at $0$) fixed
point of the SNT. It remains to check that $\mu \in \mathcal{P}_{h,m}^{+}$.
Clearly,
\begin{equation}
\underset{s\rightarrow +0}{\lim \sup }(-\varphi ^{\prime }(s))\leq m\text{.}
\end{equation}
\newline
So we should only study the lower limit.

To this end for $n=0,1,...$ put $\Phi _{n}(s):=\log (-\varphi _{n}^{\prime
}(e^{-s}))$, $\Psi _{n}(s):=\log \varphi _{n}(e^{-s}).$ Note that in view of
assumptions $\pi (dz):=-\lambda zh^{\leftarrow }(dz)$ is a probability
measure and let $\theta $, $\theta _{1}$, $\theta _{2}$, ... be independent
rv's with this distribution. Under these notations one obtains from (15) by
change of variable
\begin{eqnarray}
\Phi _{n}(s) &=&\Psi _{n}(s)+\log \int_{0}^{\infty }-\varphi _{n-1}^{\prime
}(e^{-s+z})\pi (de^{z})=  \notag \\
&=&\Psi _{n}(s)+\log \mathbb{E}\{\exp \Phi _{n-1}(s-\log \theta )\}\geq
\notag \\
&\geq &\Psi _{n}(s)+\mathbb{E}\Phi _{n-1}(s-\log \theta )\geq -me^{-s}+%
\mathbb{E}\Phi _{n-1}(s-\log \theta )\text{.}
\end{eqnarray}
Above the first inequality follows by Jensen's inequality and the second one
follows by monotonicity of $\{\Psi _{n}\}.$ \newline
Consider the random walk $S_{0}=0$, $S_{n}=-\sum_{i=1}^{n}\log \theta _{i}$,
$n=1,2,...$ On iterating (17) one gets
\begin{equation*}
\Phi _{n}(s)\geq \mathbb{E}\Phi _{0}(s+S_{n})-me^{-s}(1+\mathbb{E}%
\sum_{i=1}^{n-1}\theta _{1}\theta _{2}...\theta _{i})\text{.}
\end{equation*}
\newline
Since $\mathbb{E}\log \theta _{i}=\lambda \int_{0}^{\infty }h(z)\log
h(z)dz<0 $ then by the strong law of large numbers $S_{n}\rightarrow +\infty
$ as $n\rightarrow \infty $ a.s. Consequently by dominated convergence
\begin{equation}
\underset{n\rightarrow \infty }{\lim }\mathbb{E}\Phi _{0}(s+S_{n})=\log m.
\end{equation}
\newline
Since $\pi $ is concentrated on $(0,1)$ then $\underset{n\rightarrow \infty
}{\lim }(1+\mathbb{E}\sum_{i=1}^{n-1}\theta _{1}\theta _{2}...\theta
_{i})=(1-\mathbb{E}\theta )^{-1}$. Therefore by using (18)
\begin{equation*}
\underset{s\rightarrow +0}{\lim \inf }(-\varphi ^{\prime }(s))=\exp \{%
\underset{s\rightarrow +\infty }{\lim \inf }\underset{n\rightarrow \infty }{%
\lim }\Phi _{n}(s)\}\geq m\text{.}
\end{equation*}
\newline
This together with (16) show that $\mu ^{\ast }:=\mu \in \mathcal{P}%
_{h,m}^{+}$.

To prove uniqueness let us assume on the contrary that there exists \emph{%
another} LST $\widetilde{\varphi }(s)$ with $\underset{s\rightarrow +0}{\lim
}s^{-1}(1-\widetilde{\varphi }(s))=m$ that satisfies (4). As in Athreya
(1969, Theorem 1), set $M(s)=\dfrac{\left| \widetilde{\varphi }(s)-\varphi
(s)\right| }{s}$ for $s>0$ and obtain from (4):
\begin{equation}
M(s)\leq \int_{0}^{1}M(sz)\pi (dz)\leq ...\leq \mathbb{E}M(s\theta
_{1}...\theta _{n}).
\end{equation}
Further for any $s>0$ $M(s)\leq \left| m-s^{-1}(\widetilde{\varphi }%
(s)-1)\right| +\left| s^{-1}(1-\varphi (s))-m\right| $ which gives $%
\underset{s\rightarrow +0}{\lim }M(s)=0$. \newline
By the strong law of large numbers and bounded convergence in (19) we
conclude that $M(s)=0$ for $s>0$. It remains to recall that $\widetilde{%
\varphi }(0)=\varphi (0)=1$ which yields $\widetilde{\varphi }(s)=\varphi
(s) $. This completes the proof.

While our third auxiliary assertion is the key ingredient to the proof of
all assertions of Section 2, and in essense makes clear the connection
between fixed points of the SNT's and perpetuities of special kind (12), the
fourth one is quite simple and again can be read from (4) with some
additional explanations. In Lemma 3.3 all random variables and distributions
involved were described just above (12).

\begin{lemma}
Let for given $\nu $ as in Proposition 2.3 a r.v. $\eta $ satisfies (12).
Then $\mu $ is a fixed point of the SNT $\mathbb{T}_{h,1}$ with $%
h^{\leftarrow }(x)=\int_{x}^{b}z^{-1}\nu (dz)$, $x\in (0,b)$ and hence $h$
is subject to Condition A and $\int_{0}^{\infty }h(z)dz=1$.\newline
Conversely, if $\mu ^{\ast }$ is a fixed point of the SNT $\mathbb{T}%
_{h,\lambda }$ with
\begin{equation}
\lambda \int_{0}^{\infty }h(z)dz=1\text{ and }h(+0)=b\in (0,1]
\end{equation}
then the r.v. $\eta $ with $\mathcal{L}(\eta )=\mu ^{\ast }$ satisfies (12)
with a r.v. $A$ whose distribution $\nu $ is concentrated on $(0,b]$ and
defined as follows: $\nu (dx)=-\lambda xh^{\leftarrow }(dx)$.
\end{lemma}

\textbf{Proof.} Let us first note that if $\nu =\delta _{1}$ then $\mu
=\delta _{0}$, the case excluded by us. By the same reasoning we remove the
indicator function from the class of possible response functions in
Condition A.\newline
Suppose that the SNT $\mathbb{T}_{h,\lambda }$ has a fixed point $\mu ^{\ast
}$ and hence $\varphi ^{\ast }(s)=\int_{0}^{\infty }e^{-sx}\mu ^{\ast }(dx)$
satisfies (4), that is, $\varphi ^{\ast }(s)=$%
\begin{equation*}
=\exp \{-\lambda \int_{0}^{\infty }(1-\varphi ^{\ast }(sh(u)))du\}=\exp
\{\lambda \int_{0}^{b}(1-\varphi ^{\ast }(sz))h^{\leftarrow }(dz)\}\text{.}
\end{equation*}
In view of Lemma 3.2 condition (20) implies $m:=\int_{0}^{\infty }x\mu
^{\ast }(dx)<\infty $. Without loss of generality we may and do assume $m=1$
and therefore\newline
$\underset{s\rightarrow +0}{\lim }s^{-1}(1-\varphi ^{\ast }(s))=1$. \newline
Suppose that a r.v. $\eta $ with $\mathbb{E}\eta =1$ satisfies (12). Then
the LT $\varphi (s)=\mathbb{E}e^{-s\eta \text{ }}$solves
\begin{equation*}
\varphi ^{^{\prime }}(s)=\varphi (s)\int_{0}^{\infty }\varphi ^{^{\prime
}}(sz)\nu (dz)\text{.}
\end{equation*}
\newline
Note that $-\varphi ^{^{\prime }}(s)$ is the LST of probability measure $%
\overline{\mu }(dx)=x\mu (dx)$. By using Fubini's Theorem one has $\ln
\varphi (s)=\int_{0}^{s}[\ln \varphi (u)]^{^{\prime
}}du=\int_{0}^{s}\int_{0}^{\infty }\varphi ^{^{\prime }}(uz)\nu
(dz)du=\int_{0}^{\infty }z^{-1}\nu (dz)(\varphi (sz)-1)$ or equivalently
\begin{equation*}
\varphi (s)=\exp \{-\int_{0}^{\infty }(1-\varphi (sz))z^{-1}\nu (dz)\}.
\end{equation*}
\newline
Put $\nu (dz)=-\lambda zh^{\leftarrow }(dz)$ and note that this implies that
the statements ''$\nu $ is a probability measure on $[0,b]$'' and (20) are
equivalent. We want to verify that $\varphi ^{\ast }(s)=\varphi (s)$.
Luckily, the way of doing so mimics that of the proof of the previous Lemma
(beginning with ''As in Athreya...''), the only difference being $M(s)=%
\dfrac{\left| \varphi (s)-\varphi ^{\ast }(s)\right| }{s}$. The proof is
completed. \newline

\begin{lemma}
Assume that $\lambda \int_{0}^{\infty }h(z)dz=1$. Then for any $\alpha \in
(0,1)$ the SNT $\mathbb{T}_{h^{1/\alpha },\lambda }$ has a fixed point $\mu
_{\alpha }^{\ast }$ whose tail is given by
\begin{equation}
\mu _{\alpha }^{\ast }(x,\infty )=\int_{0}^{\infty }s_{\alpha
}(xt^{-1/\alpha },\infty )\mu ^{\ast }(dx)
\end{equation}
where $\mu ^{\ast }$ is a fixed point of $\mathbb{T}_{h,\lambda }$ with
finite mean; $s_{\alpha }$ is a strictly stable positive distribution with
index of stability $\alpha $ or equivalently
\begin{equation}
\int_{0}^{\infty }e^{-sx}\mu _{\alpha }^{\ast }(dx)=\varphi ^{\ast
}(s^{\alpha })\text{,}
\end{equation}
where $\varphi ^{\ast }(s)$ is the LST of $\mu ^{\ast }$.
\end{lemma}

\textbf{Proof.} Set $\varphi _{\alpha }^{\ast }(s)=\varphi ^{\ast
}(s^{\alpha })$ and let $m$ be the mean of $\mu ^{\ast }$. A formal
substitution in (4) $s^{\alpha }$ instead of $s$ gives
\begin{equation}
\varphi _{\alpha }^{\ast }(s)=\exp \{-\lambda \int_{0}^{\infty }(1-\varphi
_{\alpha }^{\ast }(sh^{1/\alpha }(u)))du\}
\end{equation}
\newline
which implies (22) provided the SNT $\mathbb{T}_{h^{1/\alpha },\lambda }$ is
well-defined or equivalently the integral in (23) converges for small $s$.
However the latter is easy since $\underset{s\rightarrow +0}{\lim }%
s^{-\alpha }(1-\varphi _{\alpha }^{\ast }(s))=m$ implies for some $%
\varepsilon >0$ and $s_{0}=s_{0}(\varepsilon )>0$ $\int_{0}^{\infty
}(1-\varphi _{\alpha }^{\ast }(sh^{1/\alpha }(u)))du\leq (m+\varepsilon
)s^{\alpha }$ for all $s\in (0,s_{0})$. To see that (22) is tantamount to
(21), recall that if $\varphi (s)=\mathbb{E}e^{-s\vartheta }$ than $\varphi
(s^{\alpha })=\mathbb{E}e^{-S_{\alpha }\vartheta ^{1/\alpha }}$, where $%
S_{\alpha }$ is a positive strictly $\alpha -$stable r.v.

\textbf{Proof} of Proposition 2.1 (a). Let $\gamma (a,b)$ be a r.v. with
gamma distribution with parameters $a,b>0$, that is, its probability density
function (p.d.f) is $p_{a,b}(x)=\dfrac{b^{a}}{\Gamma (a)}x^{a-1}e^{-bx}$, $%
x>0$, and $\beta (c,d)$ be a r.v. with beta distribution of the first kind
with parameters $c,d>0$, that is with p.d.f. $q_{c,d}(x)=\dfrac{\Gamma (c+d)%
}{\Gamma (c)\Gamma (d)}x^{c-1}(1-x)^{d-1}$, $x\in (0,1)$. The well-known
result due to Stuart (1962) asserts that for any positive $\alpha _{1}$, $%
\alpha _{2}$ $\gamma (1,\alpha _{2})\overset{d}{=}\beta (1,\alpha
_{1})\gamma (1+a_{1},\alpha _{2})$. This together with the obvious equality $%
\gamma (1+\alpha _{1},\alpha _{2})\overset{d}{=}\gamma (1,\alpha
_{2})+\gamma (a_{1},\alpha _{2})$ imply for $\alpha _{1}=$ $\alpha
_{2}=\alpha $
\begin{equation}
\gamma (1+\alpha ,\alpha )\overset{d}{=}\gamma (a,\alpha )+\beta (1,\alpha
)\gamma (1+\alpha ,\alpha )\text{.}
\end{equation}
\newline
It remains to note that $\mathbb{E}\gamma (a,\alpha )=1$ and hence (12) is
nothing more than (24) with $\overline{\eta }\overset{d}{=}\gamma (1+\alpha
,\alpha )$, $\eta \overset{d}{=}\gamma (\alpha ,\alpha )$ and $A\overset{d}{=%
}\beta (1,\alpha )$. By Lemma 3.3 $\gamma (a,\alpha )$ is a fixed point of
the SNT $\mathbb{T}_{h,1\text{ }}$with $h$ defined by its inverse $%
h^{\leftarrow }(x)=\int_{x}^{1}z^{-1}q_{1,\alpha }(z))dz$. To complete the
study of gamma distributions it suffices to note that fixed points of the
SNT (2) are scale invariant by Lemma 3.1.

Since $\int_{0}^{\infty }e^{-zx}\mu _{\alpha ,\beta ,\gamma
}(dx)=\int_{0}^{\infty }e^{-z^{\gamma }x}\mu _{\alpha ,\beta }(dx)$, an
appeal to Lemma 3.4 finishes the proof.\newline
b) Pitman, Yor (2000, Proposition 12(i,iii)) proved that $S2$ distribution $%
\mu _{2}$ given via the LST $\varphi _{2}(s)=\left( \dfrac{\sqrt{2s}}{\sinh
\sqrt{2s}}\right) ^{2}$satisfies (12) with $\nu =\mathcal{L}(A)$ such that $%
\nu (dx)=(x^{-1/2}-1)dx$, $x\in (0,1)$. Hence, by Lemma 3.3 $\mu _{2}$ is a
fixed point of $\mathbb{T}_{h,1}$with $h$ being defined via its inverse $%
h^{\leftarrow }(x)=\int_{x}^{1}z^{-1}\nu (dz)=\ln x+2x^{-1/2}-2$,\ $x\in
(0,1)$. An appeal to Lemma 3.1 proves that the same is true for $\mu
_{\delta }$.\newline
The conclusion regarding distributions given by (11) comes from Lemma 3.4.
This finishes the proof. \

\begin{remark}
\strut Formula (24) is well-known and especially often mentioned in the
literature on perpetuities. There are some its extensions which can be found
in Dufresne (1995, 1998).\ \newline
\end{remark}

\textbf{Proof }of Proposition 2.2. We provide an explicit example of such a
possibility. In fact, we intend to show that 1) the response function $h(u)=%
\dfrac{1}{(\cosh u)^{2}}$ generate SD fixed points of the SNT\ $\mathbb{T}%
_{h,1\text{ }}$being $\gamma (1/2,1/2)$ distributions; 2) so defined $h$ is
log-concave. \newline
To this end let us turn to Proposition 2.1 to obtain that $\gamma (1/2,1/2)$
distribution is a fixed point of $\mathbb{T}_{h,1\text{ }}$where $h$ is
given via its inverse $h^{\leftarrow
}(u)=2^{-1}\int_{u}^{1}z^{-1}(1-z)^{-1/2}dz=-\dfrac{1}{2}\ln \dfrac{%
1-(1-u)^{1/2}}{1+(1+u)^{1/2}}$, $u>0$. Now it is easily seen that the
corresponding $h$ is of the form stated above by appealing at final stage to
the well-known relation
\begin{equation*}
1-(\tanh u)^{2}=\dfrac{1}{(\cosh u)^{2}}\text{.}
\end{equation*}
\newline
Log-concavity of $h$ follows from the relation $(\ln h(u))^{^{\prime \prime
}}=\dfrac{-2}{(\cosh u)^{2}}<0$. This completes the proof.

\textbf{Proof} of Proposition 2.3. e) is our Lemma 3.3. All the other parts
of the Proposition can be obtained by appealing to e) as follows: a) is a
consequence of Lemma 3.2; b) finiteness of some exponential moments is a
part of Theorem 1.1(b) of IJ(2002); c) is quite trivial and can be read from
(5); d) is a part of Theorem 1.2 of IJ(2002).

\end{document}